\documentclass[12pt]{amsart}

\usepackage{amsfonts}
\usepackage{amscd}
\usepackage{epsfig}
\usepackage{graphicx}
\usepackage{amsmath,amssymb}
\usepackage{amsthm}
\usepackage{url}
\usepackage{color}

\newtheorem{thm}{Theorem}[section]

\newtheorem{prop}[thm]{Proposition}
\newtheorem{cor}[thm]{Corollary}

\newtheorem{rmk}[thm]{Remark}

\newtheorem{ques}[thm]{Question}
\newtheorem{ex}[thm]{Example}

\def\O{{\mathcal O}}

\def\I{{\mathcal I}}

\def\Pthree{{\mathbb P}^3}
\def\surj{\rightarrow \hskip -10pt \rightarrow}

\def\Aone{{\mathbb A}^1}

\def\Hilb{\mathop{\rm Hilb}}

\def\Hom{\mathop{\rm Hom}}

\def\Ext{\mathop{\rm Ext}}

\def\ann{\mathop{\rm Ann}}
\def\length{\mathop{\rm length}}
\def\deg{\mathop{\rm deg}}

\def\cok{\mathop{\rm Coker}}
\def\ker{\mathop{\rm Ker}}

\title{Detaching embedded points}

\author{Dawei Chen}

\address{University of Illinois at
  Chicago, Department of Mathematics, Statistics and Computer Science, Chicago, IL 60607}

\email{dwchen@math.uic.edu}

\author{Scott Nollet}

\address{Texas Christian University, Department of Mathematics, 
Fort Worth, TX 76129}

\email{s.nollet@tcu.edu}


\begin{document}
\bibliographystyle{plain}

\begin{abstract}
Suppose that $X \subset Y \subset \mathbb P^N$ differ at finitely many 
points: when is $Y$ a flat limit of $X$ union isolated points? 
Our main result says that this is possible when $X$ is a local complete 
intersection of codimension 2 and the multiplicities are at most 3. 
We show by example that no hypothesis can be weakened: the conclusion 
fails for (a) embedded points of multiplicity $> 3$, (b) local complete intersections 
$X$ of codimension $> 2$ and (c) non-local complete intersections of codimension 2. 
As applications, we determine the irreducible components of Hilbert schemes of 
space curves with high arithmetic genus and show the smoothness of the Hilbert 
component whose general member is a plane curve union a point in $\mathbb P^3$.
\end{abstract}

\maketitle

\section{Introduction}

An attractive aspect of algebraic geometry is that moduli spaces for its objects 
tend themselves to be algebraic varieties.
Ever since Grothendieck proved their existence \cite{G}, the Hilbert schemes 
$\Hilb^{p(z)} (\mathbb P^N)$ classifying flat families of subschemes in 
$\mathbb P^N$ with fixed Hilbert polynomial $p(z)$ have drawn great interest. 
One of the first major results was the connectedness of Hilbert schemes, 
proved in Hartshorne's thesis \cite{H1}. More recently Liaison theory 
\cite{MDP2,LT,PeS} has focused attention on the Hilbert schemes $H_{d,g}$ of 
locally Cohen-Macaulay curves in $\mathbb P^3$ of degree $d$ and arithmetic 
genus $g$, whose connectedness remains an open question \cite{Hartconn,N,N2,NS}. 

While Hilbert schemes can be quite complicated in general, Piene and Schlessinger 
gave a satisfying picture of $\Hilb^{3z+1} (\mathbb P^3)$: 
there are two smooth irreducible components of dimensions 12 and 15 which meet 
transversely along an 11-dimensional family \cite{PS}. The first named author 
applied the Mori program to the 12-dimensional component of twisted cubic curves, 
working out the the effective cone decomposition and the corresponding 
models, exhibiting it as a flip of the Kontsevich moduli space of stable maps 
over the Chow variety \cite{C}. Similarly the Hilbert scheme component of unions 
of a pair of codimension two linear subspaces of $\mathbb P^N$ is a smooth Mori 
dream space \cite{CCN}. 

We would like to have a similar understanding of the birational geometry of the Hilbert scheme 
component of rational quartic curves in $\mathbb P^3$. The first obstacle is a determination 
of the irreducible components of $\Hilb^{4z+1} (\mathbb P^3)$, which requires ruling out 
components whose general 
member has embedded points: this is a serious concern, for we show that such components 
exist (Example \ref{embedcomp}). This problem led us to ask the following:
\begin{ques}\label{Q1}{\em 
If $X \subset Y \subset \mathbb P^N$ differ at finitely many points, when is $Y$ a flat 
limit of $X$ union isolated points?
\em}\end{ques}

Question \ref{Q1} is already interesting when $X$ is empty. Fogarty observed that 
$\Hilb^{d}(\mathbb P^2)$ is irreducible \cite{F}, but Iarrobino showed that the Hilbert 
scheme $\Hilb^{d}(\mathbb P^3)$ is reducible for $d \gg 0$ \cite{I}. The minimum 
such value of $d$ is still unknown. Work of Cartwright, Erman, Velasco and Viray \cite{E}
extends work of Mazzola \cite{M} to prove that for $d \leq 8$, $\Hilb^{d}(\mathbb P^N)$ 
is reducible if and only if $d=8$ and $N \geq 4$. We are interested in the case 
$\dim X > 0$, when Question \ref{Q1} amounts to determining when embedded points of $Y$ 
can be detached. 

The following criterion tells when {\it all} subschemes obtained from $X$ by adding a 
single embedded point at $p \in X$ are flat limits of $X$ union an isolated point 
(see Corollary \ref{local1}). 

\begin{thm}\label{one}
For $p \in X \subset \mathbb P^N$, the following are equivalent:
\begin{enumerate}
\item All subschemes $Y$ obtained from $X$ by adding an embedded point at $p$ of 
multiplicity one are flat limits of $X$ union an isolated point.
\item The ideal sheaf $\I_X$ has $r$ minimal generators at $p$ with $r \leq N$ and 
$\pi^{-1}(p) \cong \mathbb P^{r-1}$, where 
$\pi:{\widetilde {\mathbb P^N}} \to \mathbb P^N$ is the blow-up at $X$.
\end{enumerate}
In particular, if $X$ is a local complete intersection, then any embedded point of 
multiplicity one can be detached from $X$. 
\end{thm}

One consequence is that many schemes $Y$ obtained from $X$ by adding a single embedded 
point can not be flat limits of $X$ union an isolated point. In fact, sometimes a single 
embedded point cannot be detached even if $X$ is allowed to move in the deformation 
(see Example \ref{1} below). Our main result gives conditions under which embedded points 
of various multiplicities can be removed:

\begin{thm}\label{main}
Let $X \subset \mathbb P^N$ be a local complete intersection of codimension 2. 
If $Y$ is obtained from $X$ by adding embedded points of multiplicity 
$\leq 3$, then $Y$ is a flat limit of $X$ union isolated points. 
\end{thm}

The hypotheses may seem restrictive, but Theorem \ref{main} is sharp in all aspects, 
as the following examples show. 

\begin{ex}\label{1}{\em
For $g \leq -15$, the Hilbert scheme $\Hilb^{4z+1-g} (\mathbb P^3)$ of curves with degree 4 and 
genus $g$ has an irreducible component $H$ of dimension $9-2g$ whose general member is the union 
of a multiplicity 4-line containing the triple line of generic embedding dimension 3 and an 
embedded point of multiplicity one. Details are given in Example \ref{embedcomp}.
\em}\end{ex}

\begin{ex}\label{2}{\em
There are local complete intersections $X \subset \mathbb P^N$ of codimension $> 2$ and $Y$ obtained 
from $X$ by adding an embedded point of multiplicity 2 which are not limits of $X$ union two isolated 
points. Let $X$ be the non-reduced curve in $\mathbb P^4$ with ideal $I_X = (x^2,y^2,z^2)$. 
The family of double point structures on $X$ has dimension 8, the same as the dimension of the family 
consisting of $X$ union 2 isolated points, hence the former cannot lie in the closure of the latter. 
See Example \ref{codim3} for details.
\em}\end{ex}

\begin{ex}\label{3}{\em
There are local complete intersections $X \subset \mathbb P^N$ of 
codimension 2 and $Y$ obtained from $X$ by adding an embedded point 
of multiplicity 4 which are not limits of $X$ union four isolated points. 
For fixed $X$ given by $I_X = (x^2,y^2)$ in $\mathbb P^N$, we give a 
family of such subschemes $Y$ having dimension $5N-6$, hence the general 
member cannot be a limit 
of $X$ union 4 isolated points for $N > 5$. See Example \ref{bowtie} for details.
\em}\end{ex}

\begin{rmk}\label{unique}{\em
(a) Note that for $Y$ and $X$ as in Theorem \ref{main}, there is an exact sequence 
\begin{equation}\label{standard}
0 \to \I_Y \to \I_X \stackrel{\varphi}{\to} K \to 0
\end{equation}
where $K$ is a sheaf of finite length. Moreover, it is clear from the sequence that the sheaf $K$ 
is uniquely determined by $Y$ (it is the quotient $\I_X/\I_Y$) and that two surjections 
$\varphi, \varphi^\prime$ yield the same subscheme $Y$ if and only if there exists an 
automorphism $\sigma$ of $K$ such that $\varphi^\prime = \sigma \circ \varphi$. 
The technique of proof deforms the pair $(\varphi, K)$. 

(b) We are {\it not} claiming that the embedded points can be pulled away one at a time, for
this is false by Example \ref{two}. 

(c) If $X$ is a {\it hypersurface} and $Y$ is obtained from $X$ by adding 
embedded points of {\it any} multiplicities, then $Y$ is a limit of $X$ union 
isolated multiple points. In particular, such $Y$ is a limit of 
$X$ union isolated reduced points if the multiplicities are $\leq 7$ (Proposition \ref{codim1}).
\em}\end{rmk}

Applying Theorem \ref{main} to plane curves, we deduce the following: 

\begin{cor}
For $d \geq 6$ and $\frac{(d-1)(d-2)}{2}-3 \leq g \leq \frac{(d-1)(d-2)}{2}$, 
the Hilbert scheme $\Hilb^{dz+1-g} (\mathbb P^3)$ is irreducible.
\end{cor}

In section 3 we give many other applications to space curves of low degree. 
Letting  $g=\frac{1}{2}(d-1)(d-2)$ be the genus of a degree-$d$ plane curve, we 
give the following smoothness result:

\begin{thm}
Let $H_d \subset \Hilb^{dz+2-g} (\mathbb P^3)$ be the closure of the family of plane curves 
and an isolated point. Then $H_d$ is smooth for all $d \geq 1$, hence 
isomorphic to the blow up of $\Hilb^{dz+1-g} (\mathbb P^3) \times \mathbb P^3$ 
along the incidence correspondence. 
\end{thm}

\begin{rmk}{\em
Similarly the Hilbert scheme of a hypersurface in $\mathbb P^N$ union an isolated point is 
smooth (Theorem \ref{hypersurface}), but the Hilbert scheme is not smooth at plane curves 
union certain double embedded points (Remark \ref{2pointsing}). 
\em}\end{rmk}

Regarding organization, we deal with the question of detaching an embedded point of muliplicity 
one in section 2, and with embedded points of multiplicities two or three in section 3. We give our 
applications to Hilbert schemes in section 4. We work over an algebraically closed field $k$ of 
arbitrary characteristic. 

{\bf Acknowledgements.} We thank Izzet Coskun, Daniel Erman, Robin Hartshorne and Michael Stillman 
for useful conversations. Part of the work was modified during the workshop ``Components of Hilbert Schemes'' 
at the American Institute of Mathematics, July 2010. We thank AIM and the organizers for invitation and hospitality. 

\section{Detaching a single embedded point}

In this section we study embedded point structures of multiplicity 
one on a local complete intersection $X \subset \mathbb P^N$ of codimension two. 
We also give a global result for ACM subschemes with 3-generated total 
ideal (Proposition \ref{cilimit}). 
We begin by determining when a single embedded point can be detached from an arbitrary subscheme 
$X \subset \mathbb P^N$. 

\begin{prop}\label{map}
For a proper subscheme $X \subset \mathbb P^N$, let $V \subset \Hilb^{p(z)+1}(\mathbb P^N)$ 
be the closed subset obtained from $X$ by adding a point $p$ (embedded or isolated). 
Then there is a diagram
\begin{equation}\label{diagram1}
\begin{array}{ccc}
{\widetilde {\mathbb P}}^N(X) & \stackrel{f}{\hookrightarrow} &  V \\
\downarrow \pi & & \downarrow h \\
\mathbb P^N & = & \mathbb P^N
\end{array}
\end{equation}
in which $\pi$ is the blow-up at $X$, $h$ associates a subscheme in $V$ 
to the added point, and $f$ extends the map 
$\mathbb P^N - X \to V$ given by $p \mapsto X \cup p$. 
Moreover, $f$ is injective.
\end{prop}
\begin{proof}
There is a uniform bound for the Castelnuovo-Mumford regularity 
of every ideal sheaf defining a closed subscheme with Hilbert polynomial $p(z)$, hence 
$h^{0} (\I_{Y} (m))$ is independent of $Y \in \Hilb^{p(z)+1} (\mathbb P^N)$ for sufficiently 
large $m$ and the map 
\[
Y \mapsto (H^0 (\mathbb P^N, \I_Y (m)) \subset H^0 (\mathbb P^N, \O_{\mathbb P^N} (m)))
\] 
yields a closed immersion $F: \Hilb^{p(z)+1}(\mathbb P^N) \hookrightarrow \mathbb G$ 
to a suitable Grassmann variety \cite{HM}.
Since $H^0 (\mathbb P^N, \I_Y (m)) \subset H^0 (\mathbb P^N, \I_X (m))$ has codimension one for 
$Y \in V$, the image $F(V)$ is contained in $\mathbb P (H^{0} (\mathbb P^N, \I_{X} (m)))^{\vee} \subset \mathbb G$. 
On the other hand, a standard construction \cite{PeS}[Prop. 4.1] yields a closed immersion 
${\widetilde {\mathbb P}}^N (X) \stackrel{j}{\hookrightarrow} \mathbb P (H^{0} (\mathbb P^N, \I_{X} (m)))^{\vee}$ 
and for each $p \in \mathbb P^{N}-X$ we have $j(\pi^{-1}(p))=F(h^{-1}(p))$. 
Since $V$ is closed, we obtain an injective map ${\widetilde {\mathbb P}}^N (X) \hookrightarrow V$ and 
accompanying diagram (\ref{diagram1}).
\end{proof}

\begin{prop}\label{onepoint}
In Diagram (\ref{diagram1}) above, the following are equivalent:
\begin{enumerate}
\item[(a)] $V$ is irreducible. 
\item[(b)] For each $p \in X$, $\dim_{k(p)} \I_{X,p} \otimes k(p) = r \leq N$ and 
$\pi^{-1}(p) \cong \mathbb P^{r-1}$.
\item[(c)] The map ${\widetilde {\mathbb P}}^N (X) \stackrel{f}{\to} V$ is bijective.
\end{enumerate}
\end{prop}

\begin{proof}
For each $Y \in V$, there is an exact sequence 
\[
0 \to \I_Y \to \I_X \to K_p \to 0
\]
where $K_p \cong \O_p$ is the skyscraper sheaf of length 1 supported at $p$.
For fixed $p$, the set of all such $Y$ is given by surjections 
\begin{equation}\label{point}
\phi \in \Hom(\I_X, \O_p) \cong \Hom(\I_{X,p}, k(p)) \cong \Hom(\I_X \otimes k(p), k(p))
\end{equation} 
modulo scalar. In view of Nakayama's lemma, we see that $h^{-1}(p) \cong \mathbb P^{r-1}_{k(p)}$, 
where $r$ is the minimal number of generators for $\I_X$ at $p$. 

In view of Proposition \ref{map}, the equivalence $(a) \iff (c)$ is clear. 
Condition (c) implies that $\pi^{-1}(p) \cong h^{-1}(p) \cong \mathbb P^{r-1}$ for each $p \in X$ 
and $r \leq N$ because $\pi^{-1}(p) \subset {\widetilde {\mathbb P}}^N (X)$ is a proper subset, 
which proves condition (b). Conversely if (b) holds, then for $p \in X$, we have injections 
$\mathbb P^{r-1} \cong \pi^{-1}(p) \hookrightarrow h^{-1}(p) \cong \mathbb P^{r-1}$ 
which must be surjective by reason of dimension, hence $f:{\widetilde {\mathbb P}}^N (X) \to V$ 
is bijective on the fibres over $\mathbb P^N$ and is therefore bijective. 
\end{proof}

The following result follows from the argument above, but is a stronger statement because it allows 
one to tell when all embedded point structures supported at a fixed point $p$ can be detached. 

\begin{thm}\label{local1}
For $p \in X \subset \mathbb P^N$, the following are equivalent:
\begin{enumerate}
\item Every subscheme $Y$ obtained from $X$ by adding a single embedded point at $p$ is a flat limit of 
$X$ union an isolated point.
\item $X$ satisfies condition (b) of Proposition \ref{onepoint} at $p$. 
\end{enumerate}
In particular, these conditions hold if $X$ is a local complete intersection.
\end{thm}
 
\begin{proof}
In the setting of Proposition \ref{map}, let $U = h^{-1}(\mathbb P^N - X) \subset V$ correspond to 
the subschemes obtained from $X$ by adding an isolated point. 
Note that $f({\widetilde {\mathbb P}}^N(X)) = {\overline U} \subset V$, since it is a closed subset with 
dense open subset $U$, hence for fixed $p \in X$ we have an inclusion 
$f(\pi^{-1}(p)) \subset h^{-1}(p) \cong \mathbb P^{r-1}$. 
Now condition (b) is equivalent to $\pi^{-1}(p) \cong \mathbb P^{r-1}$, which is equivalent to 
$f(\pi^{-1}(p)) = h^{-1}(p)$ by reason of dimension and irreducibility of $\mathbb P^{r-1}$, 
but this equality is equivalent to $h^{-1}(p) \subset f({\widetilde {\mathbb P}}^N (X)) = {\overline U}$, 
which is equivalent to condition (a). 
If $X$ is a local complete intersection of codimension $r$, then it is well-known that
$\pi^{-1}(p) \cong \mathbb P^{r-1}$ \cite{AG}[II, Thm. 8.24 (b)], hence condition (b) from 
Proposition \ref{onepoint} holds.
\end{proof}

We can make a stronger statement when $X$ is a hypersurface. 

\begin{prop}\label{codim1} If $Y$ is obtained from a hypersurface $X\subset \mathbb P^N$ 
by adding embedded points of any multiplicities, then $Y$ is a limit of $X$ 
union isolated multiple points. In particular, $Y$ is a limit of $X$ union isolated reduced points 
if the multiplicities are $\leq 7$. 
Moreover, for $N\geq 4$, there exist multiplicity-8 embedded structures $Y$ that cannot be limits of 
$X$ union 8 reduced points. 
\end{prop}

\begin{proof}
Suppose that $Y$ is defined by the surjection $\I_X \to K$, where $K$ is of finite 
length supported at the embedded points $p$. Now $K \cong \oplus_p \O_{Z_p}$ 
for finite length subschemes $Z_p$ supported at $p$ ($\I_X$ is principal) and locally
$\I_Y = \I_X \cdot \I_Z$. If $Z_t$ is a rigid deformation of $Z$ away from $X$ with $Z_0 = Z$, 
then $\I_{X \cup Z_t}=\I_X \cdot \I_{Z_t}$ locally for $t \neq 0$. 
From this description of the local ideals it is clear that $\displaystyle \lim_{t \to 0} X \cup Z_t = Y$.  
If the length of $Z_t$ is $\leq 7$, then $Z_t$ is a limit of reduced points \cite{E, M}, hence $Y$ 
is a limit of $X$ union isolated reduced points. 

By \cite{E}, for $N\geq 4$, there exists a non-smoothable, length-8 subscheme $Z\subset \mathbb P^N$. 
Let $Y$ be the subscheme defined by the surjection $\I_X \rightarrow \mathcal O_Z$, where the hypersurface 
$X$ contains the support of $Z$.  Since $I_X$ is generated by a single element $f$, we have $I_Y = (fI_Z)$. 
We claim that $Y$ cannot be a limit of $X$ union 8 isolated points. Otherwise $I_Y$ would be the limit of ideals 
$(fI_{Z_t})$, where $Z_t$ consists of 8 isolated points away from $X$. From the expression of $I_Y$, it would imply 
that $I_Z$ is the limit of $I_{Z_t}$. Consequently $Z$ would be a limit of 8 reduced points, contradiction. 
\end{proof}

\begin{ex}\label{nonlci}{\em
We give two examples. 

(a) If $X \subset \mathbb P^N$ is a local complete intersection of codimension $r$ at $p$, 
then $\I_{X,p} = (f_1,f_2,\dots,f_r) \subset \O_{\mathbb P^N,p}$ with $f_i$ a regular sequence 
at $p$. An embedded point is determined by a surjection $\varphi: \I_{X,p} \surj k(p)$. 
After changing generators, we may assume $\varphi(f_1)=1$ and 
$\varphi(f_i)=0$ for $i>1$ so that the local ideal for the corresponding subscheme $Y$ is 
\[
I_Y = (m_p f_1, f_2,f_3,\dots,f_r).
\]
  
(b) Let $C \subset \mathbb A^3$ be the union of the three coordinate axes, so $C$ has ideal 
$I_C = (xy,xz,yz)$. Away from the origin $C$ is a local complete intersection. 
Working on a local affine patch, one computes that the blow-up at $C$ has 
fibre $\mathbb P^2$ over the origin, so condition (b) holds at each point. 
It follows from Theorem \ref{local1} that any curve $D$ obtained from $C$ by adding an 
embedded point is a flat limit of $C$ with a varying isolated point. 
\em}\end{ex}
 
\begin{ex}\label{onepointfailure}{\em
We give two examples where Theorem \ref{local1} does not apply. 

(a) Fix a line $L \subset \mathbb P^3$. If $C$ is the curve defined by $\I_C = \I_L^d$ 
with $d > 1$, then $\I_C$ is locally generated by $d+1$ elements at $p \in C$, 
but $\pi^{-1}(p) \cong \mathbb P^1$ (the blow-up of $\mathbb P^3$ at $\I_L$ and $\I_L^d$ 
are isomorphic \cite{AG}[II, Ex. 7.11 (a)]), so condition (b) of Proposition \ref{onepoint} fails. 

(b) The planar double line with ideal $(x,y^2)$ links via the complete intersection $(x^2,y^3)$ to 
an ACM curve $C_0$ with ideal $(x^2,xy,y^3)$ and these three generators form a minimal generating 
set at each point $p \in C_0$. Here it is not possible that $\pi^{-1}(p) \cong \mathbb P^{2}$ 
for each $p \in C$, for then the exceptional divisor would have dimension 3. 
Therefore a general embedded point cannot be detached while leaving $C_0$ fixed. 
The next result deals with this problem.
\em}\end{ex}

\begin{prop}\label{cilimit}
Let $X_0 \subset \mathbb P^N$ be ACM of codimension two with minimal resolution 
\begin{equation}\label{acmres}
0 \to \oplus_{j=1}^2 \O(-b_j) \stackrel{\psi_0}{\to} \oplus_{i=1}^3 \O(-a_i) 
\stackrel{\pi}{\to} \I_{X_0} \to 0
\end{equation}
which is a flat limit of local complete intersections. 
Then each subscheme $Y$ obtained from $X_0$ by adding a point is the flat limit of local 
complete intersection ACM subschemes and a moving isolated point. 
\end{prop}

\begin{proof}
The ideal sheaf $\I_Y$ is the kernel of a surjection $\varphi: \I_{X_0} \surj \O_p$. 
The composition $\varphi \circ \pi:\oplus_{i=1}^3 \O(-a_i) \to \O_p$ lifts to 
$\displaystyle \widetilde \varphi : \oplus_{i=1}^3 \O(-a_i) \to \O$
because $H^0_* (\O_{\mathbb P^3}) \to H^0_*(\O_p)$ is surjective in positive degrees. 
Write this map as $\widetilde \varphi = (A_1,A_2,A_3) \in H^0(\oplus_{i=1}^3 \O(a_i))$. 
We claim that for the general such lift $\widetilde \varphi$, which may be written 
$(A_1+B_1,A_2+B_2,A_3+B_3)$ with $B_i \in I_p$, the composition 
$\widetilde \varphi \circ \psi_0: \oplus_{j=1}^2 \O(-b_j) \to \O$ 
vanishes along a complete intersection of two hypersurfaces. 

To see this, write 
\[
\psi_0 = \left(\begin{array}{ccc}f_1 & f_2 & f_3 \\ g_1 & g_2 & g_3 \end{array}\right),
\]
when $\widetilde \varphi \circ \psi_0$ is given by $\sum (A_i+B_i) f_i$ and $\sum (A_i+B_i) g_i$ and 
it suffices to show that the zeros of these equations meet properly. Letting $L$ be a line 
missing $X_0$ (and $p$), we will show that there are no common zeros along $L$. Restricting 
the resolution (\ref{acmres}) to $L$ and dualizing yields the exact sequence
\[
0 \to \O_L \to \oplus_{i=1}^3 \O_L(a_i) 
\stackrel{\psi_0^\vee \otimes \O_L}{\to} \oplus_{j=1}^2 \O_L(b_j) \to 0.
\]
Since $b_j > 0$, the rank two bundle on the right has a non-vanishing section, which lifts to a section  
$(r_1,r_2,r_3)$ of the rank three bundle $\oplus \O_L(a_i)$. 
Since the equations in $I_p$ of degree $d > 0$ cut out the complete linear system $H^0(\O_L (d))$, 
we can find $B_i \in I_p$ such that $(A_i+B_i)|L=r_i$ for $i=1,2,3$ and this choice proves our claim 
because the non-vanishing image of $(r_1,r_2,r_3)$ in $\oplus \O_L(b_j)$ is given by the restrictions of the 
polynomials $\sum (A_i+B_i) f_i$ and $\sum (A_i+B_i) g_i$, hence these have no common zeros along $L$.

The general map $\psi_1: \oplus_{j=1}^2 \O(-b_j) \to \oplus_{i=1}^3 \O(-a_i)$ is injective with cokernel 
the ideal sheaf of a local complete intersection ACM subscheme $X_1$. 
The linear deformation $\psi_t = (1-t) \psi_0 + t \psi_1$ yields the maps 
$\tilde \varphi \circ \psi_t: \oplus_{j=1}^2 \O(-b_j) \to \O$, which define a family 
of schemes $S_t$ which are complete intersections in a neighborhood of $t=0$ because 
this is true of $S_0$ by construction. 
If $S \subset \mathbb P^N \times \Aone$ is the total family, there is an integral 
curve $T$ through $(p,0)$ inside $S$ which is not vertical at $(p,0)$ and base extension by 
$T \to \Aone$ allows us to pick out a moving point $p_t \in S_t$ with $p_0=p$. 

Let 
$\varphi_t: \oplus_{i=1}^3 \O(-a_i) \stackrel{\tilde \varphi}{\to} \O \to \O_{p_t} = k(p_t)$ 
be the composition. 
For general $t \neq 0$, $\cok \psi_t$ is the ideal sheaf of an ACM local complete intersection $X_t$ and 
$\varphi_t \circ \psi_t = 0$ by construction 
(since $p_t \in S_t$), hence we get induced maps $\I_{X_t} \to \O_{p_t}$: since 
$\varphi_0$ is onto, so is $\varphi_t$ for general $t$, the kernels giving a family of ideals 
$\I_{Y_t}$ for a family of local complete intersections $X_t$ converging to 
$X_0$ along with a point $p_t$ converging to $p=p_0$. If $p_t \not \in X_t$, then 
we are done. If $p_t \in X_t$ for each $t$, then we have at least shown that 
$Y$ is a limit of complete intersections having an embedded point. 
By Proposition \ref{onepoint}, these are limits of local complete 
intersections with an isolated point and we again conclude.
\end{proof}

\begin{ex}\label{tripleline}{\em
The ACM curve $C_0$ with ideal $(x^2,xy,y^2)$ is the flat limit of twisted cubic curves in $\mathbb P^3$, 
hence any curve $D$ obtained from $C_0$ by adding an embedded point is a flat limit of twisted cubic curves 
and an isolated point. 
\em}\end{ex}

\begin{ex}\label{extremal}{\em
Let $C_0$ be any locally Cohen-Macaulay curve in $\Hilb^{4z}(\mathbb P^3)$. 
Then $C_0$ achieves equality for the genus bound for non-planar curves \cite{genus, MDP3}, 
hence is ACM. Since $\chi (\I_{C_0} (1)) = 0$, we have $h^2 (\I_{C_0} (1))=0$ 
and $\I_{C_0}$ is Mumford 3-regular, which implies that $\I_{C_0}(3)$ is generated by global sections. 
Also $\chi (\I_{C_0} (2))=2$ tells us that $C_0$ lies on two quadric surfaces, 
so $C_0$ is contained in a complete intersection of quadric and cubic, which links 
$C_0$ to a curve $D$ of degree $2$ and genus $0$, a plane conic. 
Using the Koszul resolution for $D$ and applying the cone construction 
\cite{LT}[Prop. 5.2.10], we obtain the resolution
\[
0 \to \O(-4) \oplus \O(-3) \to \O(-3) \oplus \O(-2)^2 \to \I_{C_0} \to 0.
\]
The family of these ACM curves is irreducible of dimension 16, the general member 
being a complete intersection of two quadrics when the ghost term splits off. 
Proposition \ref{cilimit} shows that any curve $D$ obtained from $C_0$ by adding an embedded 
point is a flat limit of complete intersections of two quadrics union an isolated point. 
\em}\end{ex}

Sometimes one cannot obtain a curve $D$ with embedded point as a limit of curves $C$ with isolated points 
even if one allows $C$ to move. In other words, the Hilbert scheme can have 
irreducible components whose general member has an embedded point:

\begin{ex}\label{embedcomp}{\em
We exhibit an irreducible component of $\Hilb^{4z+1-g}(\mathbb P^3)$ whose general member 
has an embedded point for $g \leq -15$. The irreducible components of the Hilbert schemes $H_{4,g}$ 
of locally Cohen-Macaulay curves of degree 4 and arithmetic genus $g$ are known \cite{NS}[Table III]. 
We note two typographical errors in the table, namely the family $G_5$ of double conics 
has dimension $13-2g$ instead of $13-3g$ \cite{NS}[p. 189] and the general member of family 
$G_{7,a}$ should be $W \cup_{3p} L$ instead of just $W$. 
Now consider the irreducible component $G_4$ of dimension $9-3g$, consisting of thick quadruple lines. 
Each curve $C \in G_4$ has a supporting line $L$ and there is an exact sequence 
\[
0 \to \I_C \to \I_W \stackrel{\phi}{\to} \O_L (-g-1) \to 0
\]
where $W$ is the triple line given by $\I_W = \I_L^2$ \cite{NS}[Prop. 2.1]. 
The surjection $\phi$ factors through $\I_W \otimes \O_L \cong \O_L(-2)^2$, hence is given by 
$\phi(x^2)=a,\phi(xy)=b,\phi(y^2)=c$ for three homogeneous polynomials $a,b,c$ of degree $-g+1$. 
Writing the total ideal of $C$ as 
\[
I_C = (x^3,x^2y,xy^2,y^3,axy-bx^2,by^2-cxy),
\]
we see that at general $p \in L$, $a,b,c$ become units in the local ring $\O_{\mathbb P^3,p}$, hence 
the ideal reduces to $I_C=(x^3,axy-bx^2,by^2-cxy)$, so $I_C$ is generically 3-generated for 
general $\phi$. 

Now consider the closed subscheme $V \subset \Hilb^{4z+2-g}(\mathbb P^3)$ 
obtained by adding an isolated or embedded point to $C$ as above, as in Proposition \ref{map}. 
The closure of the component corresponding to $C$ along with isolated points has dimension 3. 
Since $I_C$ is generically 3-generated, the set of embedded point structures at general 
$p \in C$ is parametrized by a $\mathbb P^2$ and we obtain a second 3-dimensional family. 
Thus $V$ is reducible with at least these two 3-dimensional components (conceivably the point 
where $I_C$ is generated by more elements could generate another family). 
Varying the curve $C \in G_4$, we obtain at least two corresponding families of dimension $12-3g$ 
(because $\dim G_4 = 9-3g)$. Let $F$ be the family whose general curve has an embedded point.

We claim that $F$ is an irreducible component of $\Hilb^{4z+2-g} (\mathbb P^3)$. The general curve $D \in F$ 
cannot be a limit of curves with having more than one isolated/embedded point (counted with multiplicity). 
Since $G_4$ is an irreducible component of $H_{4,g}$ for $g \leq -2$, $D$ is not a flat limit of 
another family of curves with a single isolated or embedded point, because the underlying locally 
Cohen-Macaulay curve $C \subset D$ is not. Finally $D$ cannot be a limit of locally Cohen-Macaulay curves 
of genus $g-1$ because the maximal dimension of such a family for $g \leq -15$ is $12-3g = \dim F$. 
\em}\end{ex}

\section{Detaching two or three points}

In this section we prove that if $Y$ has embedded points of multiplicity two (Proposition \ref{twopoint}) or three (Proposition \ref{threepoint}) 
and the underlying locally Cohen-Macaulay subscheme $X \subset \mathbb P^N$ is a local complete intersection of codimension two, 
then $Y$ is a flat limit of $X$ and isolated points. Along with Theorem \ref{local1}, this completes the proof of our main result Theorem \ref{main}, which says that
embedded points of multiplicity $\leq 3$ can be detached from a local complete intersection of codimension two. 

We begin with two propositions which take care of most of the easier cases, dealing with the more 
difficult cases in the main result of this section. We give examples to 
show that these results fail for local complete intesections of 
codimension $> 2$ (Example \ref{codim3}) and for embedded points of multiplicity 
$ \geq 4$ (Example \ref{bowtie}). 

\begin{prop}\label{embdim1}
Let $X \subset \mathbb P^N$ be a local complete intersection of codimension two, $Z$ a zero 
dimensional subscheme of embedding dimension $\leq 1$ and suppose that $Y$ is defined by the 
exact sequence 
\[
0 \to \I_Y \to \I_X \stackrel{\varphi}{\to} \O_Z \to 0
\]
Then $Y$ is the flat limit of $X$ union isolated points.
\end{prop}

\begin{proof}
Since the result is local, we may assume that $Z$ is supported at a point $p$ and has length $d$. 
Since $Z$ has embedding dimension $\leq 1$, we can choose a smooth connected curve $C_0$ of high degree containing 
$Z$ and not entirely in $X$. If $p \not \in X$, the result is clear because $Z$ is a flat limit of isolated points on $C_0$. 
In the interesting case $p \in X$, our idea is to take a deformation $C_t$ of $C_0$ and use $d$ isolated points in $C_t$ to realize the 
detaching process. 

Let $C$ be a rigid deformation of $C_0$ (induced by PGL($N+1$)) which misses $X$. 
Now for $m \gg 0$, the general pair $F,G \in H^0(\I_X (m))$ give hypersurfaces which cut out $X$ 
in an open neighborhood of $p$; for the purposes of this proof we may assume that $X$ is equal to 
the complete intersection defined by $F$ and $G$, giving the Koszul resolution

\begin{equation}\label{koszul}
0 \to \O(-2m) \stackrel{\psi}{\to} \O(-m) \oplus \O(-m) \stackrel{\pi}{\to} \I_C \to 0.
\end{equation} 

Because the restriction map $H^0(\O(m)) \to H^0(\O_Z (m))$ is surjective for $m \gg 0$, 
we can lift the images of $F,G$ to $\O=\O_{\mathbb P^N}$, hence the composition 
$\varphi \circ \pi: \O(-m)^2 \surj \O_Z$ factors through $\O$ and we obtain 
$\widetilde \varphi: \O(-m)^2 \to \O$ inducing $\varphi$. 
The composition $\widetilde \varphi \circ \psi$ vanishes on a hypersurface $S$ of degree $2m$ containing 
both $X$ and $Z$. 

By Bertini's theorem, we could have 
chosen the equations $F,G$ cutting out $X$ near $p$ to be smooth away from $X$, 
meeting $C$ in disjoint reduced sets of points, so the restrictions to $C$ 
induce a sheaf surjection $\O_{C}^2 \to \O_{C}(m)$. 
If $\widetilde \varphi$ is given by $A_0,B_0 \in H^0(\O(m))$, then $S$ has equation 
$FA_0 + GB_0 = 0$, but $A_0,B_0$ are only determined up to elements of $H^0 (\I_Z (m))$. 
Since the natural map $H^0(\I_Z (m))^2 \to H^0 (\O_{C}(m))$ is surjective, given by 
$(A,B) \mapsto FA+GB$, we may choose $A_0, B_0$ to assume that $S \cap C$ is a reduced 
set of $2m(\deg C)$ points. 

Now consider a rigid deformation $C_t$ from $C$ to $C_0$, parametrized by $t \in \Aone$. 
Now $C_0 \cap S$ contains $Z$ at $p$ and $C_t \cap S$ consists of $2m(\deg C)$ reduced points 
for general $t \neq 0$. Possibly after a base extension, we may pick $d$ distinct points 
$p_{1,t}, p_{2,t}, \dots p_{d,t}$ in $S \cap C_t$ near $p$. 
Letting $Z_t = \{p_{1,t}, p_{2,t}, \dots p_{d,t}\}$, the flat limit of $Z_t$ is exactly 
$Z$, because the ideal of the limit contains the equations of the curve $C_0$ by construction, 
and $Z$ is the unique length $d$ subscheme of $C_0$ at $p$. 

Letting $\varphi_t$ be the composition 
$\O(-m)^2 \stackrel{\widetilde \varphi}{\to} \O \to \O_{Z_t}$, 
we have $\varphi_t \circ \psi =0$ by construction, so these maps factor through 
$\I_X$ and we obtain a family of maps $\varphi_t:\I_X \to \O_{Z_t}$. 
Since $\varphi_0=\varphi$ is surjective, so are $\varphi_t$ for $t$ near $0$ 
and the family $\I_{Y_t} = \ker \varphi_t$ gives the desired family. 
\end{proof}

\begin{prop}\label{pullone}
Let $X \subset \mathbb P^N$ be a local complete intersection, 
$K$ a sheaf of finite length supported at $p$, and $Y, Y^1$ defined by the 
commutative diagram of short exact sequences
\[
\begin{array}{ccccc}
\I_Y & \to & \I_X & \stackrel{\varphi=(\alpha,\beta)}{\to} & K 
\oplus \O_p \\
 \downarrow & & \downarrow & & \downarrow \pi_1 \\
\I_{Y^1} & \to & \I_X & \stackrel{\alpha}{\to} & K  
\end{array}
\]
Then $Y$ is the flat limit of $Y^1$ union an isolated point. In particular, if $Y^1$ is the 
flat limit of $X$ and isolated points, then so is $Y$.
\end{prop}

\begin{proof}
The result is local at $p$. The direct sum allows us to write 
$\varphi = (\alpha, \beta)$, where $\alpha$ defines 
$Y^1$ as above. 
The surjection $\beta:\I_X \to \O_p$ defines an embedded point 
structure $Y^2$ on $X$. Since $X$ is a local complete intersection, $Y^2$ 
is a flat limit of $X$ union an isolated point by Theorem \ref{local1}, 
meaning that there is a flat family $Y^2_t$ for $t \in T$ with $Y^2_0=Y^2$ 
and $Y^2_t = X \cup p_t$ with $p_t \not \in X$ for $t \neq 0$. This gives 
a family of surjections $\beta_t:\I_X \to \O_{p_t}$ with 
$\I_{Y^2_t}=\ker \beta_t$ and $\beta_0=\beta$.

Let $\gamma_t:\I_{Y^1} \subset \I_X \stackrel{\beta_t}{\to} \O_{p_t}$ be 
the composite map. Clearly $\gamma_t$ is surjective for $t \neq 0$, 
because $p_t \not \in \I_X$, so the inclusion $\I_{Y^1} \subset \I_X$ is 
an equality at these points. The map $\gamma_0$ is also a surjection, 
since locally at $p$, if we choose $f \in \I_X$ which maps to $(0,1)$, 
then $\alpha(f)=0 \Rightarrow f \in \I_{Y^1}$ and $\gamma_0(f)=1$.
This family of maps gives a flat family $Y_t$ and for $t \neq 0$, $Y_t$ 
consists of $Y^1$ union an isolated point. Finally, the kernel of 
$\gamma_0:\I_{Y^1} \to \O_p$ is exactly $\I_Y$, for 
$g \in \I_{Y^1} \Rightarrow g \in \I_X$ and $\alpha(g)=0$. 
Now $\gamma_0 (g)=0 \iff \beta(g)=0 \iff \varphi(g)=0 \iff g \in \I_Y$.
\end{proof}

\begin{prop}\label{twopoint}
Let $X \subset \mathbb P^N$ be a local complete intersection of 
codimension two and obtain $Y$ by adding a double embedded point with associated exact sequence
\[
0 \to \I_Y \to \I_X \stackrel{\varphi}{\to} K_p \to 0
\]
with $K_p$ a length 2 sheaf supported at $p$. Then either 
\begin{enumerate}
\item[(a)] $K_p \cong \O_p \oplus \O_p$ or 
\item[(b)] $K_p \cong \O_Z$, where $Z \subset \mathbb P^N$ has length two.
\end{enumerate}
In either case, $Y$ is the flat limit of $X$ union two isolated points. 
\end{prop}

\begin{proof} 
If $K_p \cong \O_p \oplus \O_p$, apply Proposition \ref{pullone}. Since 
$Y^1$ is obtained from $X$ by adding a single embedded point, it is a 
limit of $X$ union an isolated point, hence $Y$ is a limit of $X$ union 
two isolated points. 

Now suppose $K_p \not \cong \O_p \oplus \O_p$. Then the surjection $K_p \to K_p \otimes k(p)$ 
is not an isomorphism, thus $K_p \otimes \O_p$ is one dimensional as an $\O_p=k(p)$ vector space. 
Therefore $K_p$ is principal by Nakayama's lemma, so there is a surjection 
$\O \to K_p$ whose kernel is the ideal sheaf $\I_Z$ of a length two subscheme, 
which is contained in a unique line and has embedding dimension one. We 
apply Proposition \ref{embdim1} to see that $Y$ is a flat limit of $X$ 
union two isolated points. 
\end{proof}

\begin{rmk}\label{mult2}{\em
We give the local equations of the embedded point structures for parts (a) 
and (b) of Proposition \ref{twopoint} for $X \subset \mathbb P^{N}$:

(a) If $I_X = (f,g)$, then $I_Y = m_p \cdot I_X$. 

(b) Replacing generators so that $\varphi (f) = 1$ and $\varphi(g)=0$, we obtain $I_Y = (g, f \cdot I_Z)$, 
$Z$ being the length two subscheme. 
\em}\end{rmk}

\begin{ex}\label{two}{\em
In case (b) of Proposition \ref{twopoint}, there is a unique subscheme $X \subset E \subset Y$ 
with an embedded point of multiplicity one, because the unique length one quotient of $\O_Z$ is 
$\O_p$, obtained by modding out by the maximal ideal. This explains why it was necessary to 
prove case (b) by pulling away two points simultaneously. 
For example, let $C \subset \mathbb A^3$ have ideal $I_C = (x^2,y^2)$ and let $p = (0,0,0)$. 
Add an embedded point to $C$ at $p$ by using the map $I_C \to k$ by 
$x^2 \mapsto 1, y^2 \mapsto 0$ to obtain $E$ with $I_E = (y^2,x^3,x^2y,x^2z)$ 
being 4-generated. 
By Proposition \ref{onepoint}, one can add a second point at $p$ to obtain $D$ with a 
double embedded point which is not limit of $E$ union an isolated point. On the other hand, 
in case (a), the structure $K_p \cong \O_p \oplus \O_p$ has many 
length 1 quotients, so it is possible to pull away the embedded points one at a time.  
\em}\end{ex}

\begin{ex}\label{codim3}{\em
Proposition \ref{twopoint} may fail for local complete intersections of 
codimension greater than two. For example, suppose that $C \subset \mathbb P^4$ is the 
complete intersection with ideal $I_C = (x^2,y^2,z^2)$. If $x,y,z,u,w$ are the projective coordinates, 
then on the affine patch $w=1$ the ideal of $C$ is also $(x^2,y^2,z^2)$. Consider the double point 
structures $D$ on $C$ given by surjections $\phi: (x^2,y^2,z^2) \to K = \O_Z$, where $Z$ is the double 
point supported at the origin with ideal $I_Z = (x,y,z,u^2)$. 
An arbitrary map $\phi:I_C \to \O_Z$ is given by 
\[
\phi(x^2)=a+bu, \;\;\; \phi(y^2) = c+du, \;\;\; \phi(z^2)=e+fu
\]
with $a,b,c,d,e,f \in k$ and any tuple $(a,b,c,d,e,f)$ is possible because $I_C \subset I_Z$. 
The automorphisms of $K=\O_Z$ are given by multiplication by $A+Bu$ with $A \neq 0$. 
Thus the maps for which $\phi(x^2)$ generates $K$ (i.e. $a \neq 0$) are 
uniquely determined up to automorphisms of $K$ by the quotients $\phi(y^2)/\phi(x^2) = (c+du)/(a+bu)$ and 
$\phi(z^2)/\phi(x^2)=(e+fu)/(a+bu)$, so these quantities uniquely determine the corresponding 
subschemes $D$ by Remark \ref{unique} (a). 
In other words, if we compose the general map above by the the automorphism of $\O_Z$ given by multiplication by 
$(a+bu)^{-1}=\frac{a-bu}{a^2}$, we may assume that $\phi(x^2)=1$ when the corresponding ideal of $D$ is given by 
\[
(x^2(I_Z),y^2-(c+du)x^2,z^2-(e+fu)x^2)=(x^3,x^2y,x^2z,x^2w^2,y^2-(c+du)x^2,z^2-(e+fu)x^2)
\] 
and each tuple $(c,d,e,f) \in k^4$ yields a distinct subscheme $D$, so we obtain a 4-dimensional 
family of such double point structures $D$. 

Finally, the same argument applies to any double point structure $D$ on $C$. Since there is a choice of 
any point $p \in C$ for the support of $K=\O_Z$ and the structure of $Z$ is uniquely determined by a 
line through $p$ (parametrized by a hyperplane $\mathbb P^3$), the family of such double point structures 
has dimension 1+3+4=8. The general such structure cannot be a limit of $C$ and two isolated points, 
for this family also has dimension 8.  
\em}\end{ex}

\begin{prop}\label{threepoint}
Let $X \subset \mathbb P^N$ be a local complete intersection of codimension two. 
Let $Y$ be the subscheme obtained from $X$ by an exact sequence
\[
0 \to \I_Y \to \I_X \stackrel{\varphi}{\to} K \to 0
\]
where $K$ is a length 3 sheaf supported at $p$. 
Then one of the following holds:
\begin{enumerate}
\item[(a)] $K \cong \O_p \oplus \O_Z$ where $Z \subset \mathbb P^N$ is a double point on a line. 
\item[(b)] $K \cong \O_Z$, where $Z \subset \mathbb P^N$ is a triple point on a line. 
\item[(c)] $K \cong \O_Z$, where $Z \subset \mathbb P^N$ is a triple point on a smooth conic. 
\item[(d)] $p \cong \O_Z$, where $Z$ is contained in a plane $H$ and 
$\I_{Z,H} = \I_{p,H}^2$.
\item[(e)] $K \cong \Hom_{\O_p}(\O_Z, \O_p)$ with $Z$ as in part (d).
\end{enumerate}
In each case, $Y$ is a flat limit of $X$ and three isolated points. 
\end{prop}

\begin{proof}
If $K$ is a direct summand, one summand is $\O_p$ and the other is 
$\O_p^2$ or $\O_Z$ for a double point $Z$: the first is not possible as a quotient 
of the locally 2-generated ideal $\I_C$, leading to case (a).
If $K$ is principal, then the surjection $\O \to K$ shows that $K \cong \O_Z$ 
for a length 3 subscheme supported at $p$ and since $h^0 (\O_Z (1))=3$ and $h^0(\O(1))=N+1$, 
$Z$ is a planar triple point: these are easy to classify, leading to cases (b)-(d). 
If $K$ is not principal, then it is 2-generated as a quotient of $\I_X$ via $\varphi$. 
The two generators have a common non-zero multiple, otherwise they would express $K$ 
as a direct sum of two principal modules. The common non-zero multiple is therefore a 
generator of the dual $\Hom(\O_Z, \O_p)$, where $\O_Z$ must be one of cases (b) - (d). However,  
the duals to cases (b) and (c) are principal and we are left with the dual of case (d), which is case (e).   

That $Y$ is a flat limit of $X$ union 3 isolated points follows from 
Proposition \ref{pullone} in case (a) and from Proposition \ref{embdim1} 
in cases (b) and (c). Cases (d) and (e) require new ideas. 

In case (d) we have $K_p \cong \O_Z$, where $Z \subset H$ is the planar triple 
point supported at $p$ of embedding dimension two. As in Proposition \ref{twopoint}, 
$X$ is contained in hypersurfaces with equations $F,G$ of degree $m \gg 0$, giving a Koszul 
resolution (\ref{koszul}), $\varphi: \I_X \to \O_Z$ lifts to 
$\tilde \varphi:\O(-m)^2 \to \O$ and there is a hypersurface $S$ of 
degree $2m$ where $\tilde \varphi \circ \psi = 0$. 
The intersection $H \cap S$ contains an integral curve $T$ passing through $p$. Our idea is to realize this triple 
embedded structure as the limit of a fixed double embedded structure at $p$ union a single one varying in $T$.

Let $\tilde T \stackrel{f}{\to} T \subset H \cong \mathbb P^2$ be the normalization of $T$ and 
choose a point $0 \in T$ such that $f(0)=p$. For $t \neq 0$, let $L_t \subset H$ be the line through $p$ and $f(t)$. 
As $t \to 0$, $f(t) \to p$ and the line $L_t$ has a unique limit $L_0$ 
(complete the assocated map $T - \{0\} \to (\mathbb P^2)^\vee$ to 
obtain this limiting line). Choose local coordinates $x,y$ on 
$\mathbb A^2 \subset \mathbb P^2$ so that $p = (0,0)$ and $L_0 = \{x=0\}$. 
The double point $W$ at $p$ with ideal $(x^2,y)$ is a closed 
subscheme of $Z$ (which has ideal $(x^2,xy,y^2)$). 
We now show that $\displaystyle \lim_{t \to 0} f(t) \cup W = Z$ 
in the Hilbert scheme of length 3 subschemes of $H$: 
if $f(t)=(a(t),b(t))$ in the local coordinates above, then the ideal 
for $W \cup f(t)$ is 
\[
I_t = (x^2,y) \cap (x - a(t), y - b(t))
\]
which contains the product of the two ideals. 
Since $\lim_{t \to 0} (a(t),b(t)) = (0,0)$, the limit 
ideal contains $(x^3,xy,y^2)$. If the line $L_t$ has equation $l_t = 0$, 
then $l_t^2 \in I_t$ and by choice of coordinates we have $\lim l_t^2 = x^2$, 
so the limit ideal also contains $x^2$ and hence $(x^2,xy,y^2)$, which defines $Z$. 

The rest is analogous to Proposition \ref{embdim1}. The composite map 
\[
\O(-2m) \stackrel{\psi}{\to} \O(-m)^2 \stackrel{\widetilde \varphi}{\to} \O 
\to \O_S \to \O_{W \cup f(t)}
\]
is zero, inducing a family of maps $\varphi_t:\I_X \to \O_{W \cup f(t)}$. 
Since $\varphi_0$ is onto, so is $\varphi_t$ for $t$ near $0$. Therefore 
the kernels $\I_{Y_t}$ give a family whose limit is 
$Y$ as $t \to 0$. Using our earlier results, for $t \neq 0$ each $Y_t$ is 
a limit of $X$ and isolated points, and therefore so is $Y$. 

Finally consider case (e), where $K_p \cong \Hom_{\O_p}(\O_Z,\O_p)$ with 
$\O_p = k(p)$ the residue field at $p$ and $Z \subset H \subset \mathbb P^N$, with 
$H$ a plane and $\I_{Z,H}=\I_{p,H}^2$. 
Choose affine coordinates $x,y,z_1,\dots,z_{N-2}$ centered at $p$ so that 
$I_H = (z_1,\dots,z_{N-2})$ and $x,y$ are coordinates for $H \cong \mathbb A^2$. 
Let $f,g$ be the restrictions of $F,G$ to this affine patch, so that $I_X = (f,g)$.
If $g - uf =h \in I_H$ for some unit $u$ in the local ring, replace $g$ with $h$ as a generator for $I_X$. In this way we 
may assume $g \in I_H$ or $(f, g)$ is not principal modulo $I_H$ locally around $p$.  
Now $\O_Z$ is generated by $1,x,y$ as an $\O_p$-vector space, so $K_p$ is generated by 
dual basis $x^*,y^*,1^*$ as a vector space and by $x^*,y^*$ as an $\O_H$-module with 
structure given by $x x^* = 1^* = y y^*, x y^* = y x^* =0$. 
Since $\varphi$ is surjective, $\varphi(f)=ax^*+by^*+c1^*, \varphi(g)=dx^*+ey^*+f1^*$
are also module generators for $K_p$, and in particular $ae-bd \neq 0$. Now consider 
the new coordinates $X = ay-bx, Y=ex-dy$ for $H$. With these one sees that 
$\ann(\varphi(f))=(I_H,X)$, $\ann(\varphi(g))=(I_H,Y)$ and $Y \varphi(f) = (ae-bd) 1^* = X \varphi(g)$.
It follows that $I_Y=(I_H(f,g), Xf, Yg, Yf - Xg)$. So by replacing the coordinates, we can present the ideal 
of $Y$ as 
\[
I_Y = (I_H(f,g), xf, yg, yf - xg).
\]

We will directly deform this ideal to obtain the result. The locus 
\[
S=\{(A,B,C,D):f(A,B)=0, g(C,D)=0, (B-D)g(C,B)-(C-A)f(C,B)=0\}
\]
contains $(0,0,0,0)$ and each component has dimension $\geq 1$, hence $S$ contains an 
integral curve $T$ through the origin: let $\sigma:T \to H \times H \cong \mathbb A^4$ be 
the inclusion with coordinate functions $\sigma(t) = (a(t),b(t),c(t),d(t))$ and 
$0 \in T$ chosen so that $\sigma(0)=(0,0,0,0)$. We claim that $T$ can be chosen with 
$(a(t),b(t)) \neq  (c(t),d(t))$. This is clear if $g \in I_H$, for then the second 
equation $g(C,D)=0$ puts no restriction on $C,D$ and $S$ is defined by only two equations: 
on a surface there are many integral curves $T$ through the origin. The other possibility by our assumption
is that $g \neq uf$ modulo $I_H$ for any invertible $u$ in an affine neighborhood of the origin. Here the restrictions of $f,g$ to $H$ have a greatest common divisor $h$ so that $f=h \cdot f_1, g=h \cdot g_1$ with $f_1,g_1$ vanishing at the origin and relatively prime modulo $I_H$ locally around the origin. If we look at the sublocus of $S$ defined as above with $f_1,g_1$ in place of $f,g$, the condition of the claim holds and we obtain the desired integral curve $T$.
  
Now consider the family of ideals 
\[
I_t = (I_H(f,g), (x-c(t)) f, (y-b(t)) g, (y-d(t)) f - (x - a(t)) g).
\]
We claim that the ideal $I_t$ scheme-theoretically cuts out exactly $X$ and the three points 
$(a(t),b(t)),(c(t),b(t)),(c(t),d(t))$ (which may be isolated or embedded, two may coincide 
if $a(t)=c(t)$ or $b(t)=d(t)$) for generic $t$ near $0$.  

The claim holds off of $H$ via the generators $I_H (f,g)$. 
At points $(x,y) \in H$ away from $(a,b),(c,b),(c,d)$ (we suppress the variable $t$) the claim also holds: 
if $x \neq c$, then $x-c$ is a unit, $f \in I_t$ and there are two cases: 
if $x = a$, then $y \neq b$, hence $y-b$ is a unit and $g \in I_t$; otherwise $x \neq a$ and 
the last equation shows that $g \in I_t$. If $x = c$, then $y \neq b,d$, so $g \in I_t$ 
and $f \in I_t$ by the last equation. 

Finally we consider $(x,y) \in \{(a,b),(c,b),(c,d)\}$. The claim is easily checked if these 
points are distinct, ($a \neq c$ and $b \neq d$) by checking that $\length I_X/I_t = 1$.
For example, at $(x,y)=(a,b)$ we have $x \neq c$ so $f \in I_t$, when $I_X/I_t$ is generated 
by $g$ alone, and since $I_H g, (y-b)g, (x-a)g \in I_t$, we have $I_X/I_t \cong k$.
The other points $(x,y)=(c,b),(c,d)$ are similar. In the degenerate case $a=c$, we need to 
show that $\length I_X/I_t = 2$ at $(x,y)=(a,b)=(c,b)$. Here $y \neq d$ so $u=(y-d)$ is a unit and 
$uf - (x-a)g \in I_t$, 
showing that $I_X/I_t$ is generated by $g$. Further $I_t$ contains $I_H g, (y-b) g$ and 
$(x-c)^2 g$ (use $(x-c)f$ and $uf - (x-c)g$), so the quotient has length two. The other degenerate case 
$b=d$ can be verified similarly. This proves the claim.

With the claim, the ideal $I_t$ cuts out $X$ and 3 other points (possibly embedded in $X$, but not all supported at the same point). Using our earlier results, these schemes are limits of $X$ and isolated points. 
Since $\displaystyle \lim_{t \to 0} (a(t),b(t),c(t),d(t))=(0,0,0,0)$ by construction, we also have 
$\displaystyle \lim_{t \to 0} I_t = I_Y$ and conclude. 
\end{proof}

\begin{rmk}\label{mult3}{\em
For $I_X=(f,g)$ near $p$ in Proposition \ref{threepoint}, we write  
local equations for the embedded point structure $Y$ according to the various 
cases: 

(a) If $K_p = \O_p \oplus \O_Z$ and $\varphi(f)=(1,0)$, $\varphi(g)=(0,1)$, then 
$I_Y = (f m_p, g I_Z)$ with $f \in I_Z$.

(b) For $K_p = \O_Z$ and $\varphi(f)=1, \varphi(g)=0$, $I_Y = (f I_Z, g)$ with $g \in I_Z$. 

(c) Similarly we have $I_Y = (f I_Z, g)$ with $g \in I_Z$. 

(d) Again we have $I_Y = (f I_Z, g)$ with $g \in I_Z$. 

(e) This is the most interesting structure. 
As shown in the proof, $I_Y = (xf-yg,yf,zf,xg,zg)$ for suitable coordinates $x,y,z$. 
\em}\end{rmk}

\begin{ex}\label{bowtie}{\em
These results show that embedded points of multiplicity at most 3 can 
always be detached from a local complete intersection of codimension 2 in $\mathbb P^N$. 
This is not possible for all embedded points of multiplicity 4. 
For linearly independent variables $X,Y,Z,W$, consider the $R=k[X,Y,Z,W]$-module given by 
\[
K = \langle a,b \rangle / (Z a, W a, X b, Y b, X a - Z b, Y a - W b).
\]
In changing the choice of vector space basis for the linear forms 
$[X,Y,Z,W]$, we obtain a family of such modules on which the group GL(4) acts. 
It's easily checked that the $R$-module automorphisms of any fixed $K$ 
have dimension 5 (one can write them down explicitly). For another
$K^\prime$ determined by basis $X^\prime, Y^\prime, Z^\prime, W^\prime$ 
and an isomorphism $\psi:K \to K^\prime$, the map $\psi$ 
uniquely determines $X^\prime, Y^\prime, Z^\prime, W^\prime$ in terms of 
$X,Y,Z,W$, because the relations yield 16 equations in 16 unknowns. One can check that 
the family of candidate isomorphisms $\psi$ has dimension 12 and a 
5-dimensional subspace corresponds to the identity coordinate change. Hence, 
we find that the isomorphism classes of such modules 
$K$ has dimension $16-(12-5)=9$. 

Now for $X \subset \mathbb P^N$ given by $I_X=(x^2,y^2)$, the family of 
embedded point structures on $X$ given by such $K$ has dimension $5N-6$:
the choice of the support of $K$ has dimension equal to $\dim X = N-2$, choosing 
the linear subspace $[x,y,z,w]$ at $p$ is given by Gr$(4,N)$ of dimension 
$4N-16$, choosing the isomorphism class of $K$ is 9 (see above), the choice 
of map $\varphi:\I_X \to K$ is 8 parameters, but the resulting family of 
ideals $\I_Y$ given by the kernels has dimension 3, because the 
automorphisms of $K$ have dimension 5. All in all, the family has 
dimension $(N-2)+(4N-16)+9+(8-5)=5N-6$. For $N \geq 6$, this is $\geq 4N$ so 
the family cannot lie in the $4N$-dimension closure of those obtained by 
unions of $X$ and isolated points.  
\em}\end{ex}

\section{Applications to Hilbert schemes}

In the previous section we proved various results about when a locally complete intersection $X$ with 
embedded points are flat limits of $X$ union isolated points. In this section 
we apply these results to describe the irreducible components of certain Hilbert 
schemes. In view of Proposition \ref{codim1}, we deduce the following:

\begin{thm}\label{hypersurface}
Let $p(z)$ be the Hilbert polynomial of a degree $d$ hypersurface in $\mathbb P^N$. Then 
\begin{enumerate}
\item[(a)] The Hilbert schemes $\Hilb^{p(z)+e} (\mathbb P^N)$ are irreducible for $0 \leq e \leq 7$. 
\item[(b)] The Hilbert scheme $\Hilb^{p(z)+1} (\mathbb P^N)$ is smooth, isomorphic to 
$\Hilb^{p(z)} (\mathbb P^N) \times \mathbb P^N$. 
\end{enumerate}
\end{thm}

\begin{proof}
It follows from Proposition \ref{codim1} that any (multiple) embedded point can be detached from a 
hypersurface and for $e \leq 7$ we also know that any subscheme of length $e \leq 7$ is a 
limit of reduced points \cite{E,M}. Therefore 
$\Hilb^{p(z)+e} (\mathbb P^N)$ is the closure of the open subset formed by a degree $d$ hypersurface and 
$e$ isolated points and $\Hilb^{p(z)+e} (\mathbb P^N)$ is irreducible of dimension ${{d+N} \choose {d}}-1+Ne$. 

Now take $e=1$. It is easily checked that the Hilbert scheme is smooth at points corresponding 
to a hypersurface and an isolated point. If $C \subset \mathbb P^N$ is a degree $d$ hypersurface and 
$D$ is obtained from $C$ by adding an embedded point located at $x_1=x_2=\dots=x_N=0$, then the 
total ideal of $D$ is simply $I_D = (x_1,x_2,\dots,x_N) \cdot I_C$ (Remark \ref{codim1}), so 
$\I_D$ is generated in degree $d+1$. 
Since the generator of $\I_C \to K$ is onto, $H^1(\I_D (n))=0$ for $n \geq d$ and so the comparison 
theorem of Piene and Schlessinger \cite{PS} applies. Now the argument \cite{PS}[Case (iii) of Lemma 4] 
goes through, which we include for self-containment: $H^0 ({\mathcal N}_D) = \Hom(I_D,A)_0$, 
where $A = S/I_D$, $S$ is the coordinate ring 
of $\mathbb P^N$, and ${\mathcal N}_D$ is the normal sheaf to $D$.
Given the dimension of $\Hilb^{p(z)+1} (\mathbb P^N)$, it suffices to prove that 
$\dim \Hom (I_D,A) \leq {{d+N} \choose {d}}-1+N$. Setting $\overline A = S/I_C$ and 
$K=I_C/I_D$, the $S$-module $K$ has Koszul resolution of the form
\[
0 \to S(-d-N) \to \dots S(-d-2)^{N(N-1)/2} \to S(-d-1)^N \to S(-d) \to K \to 0.
\]
Applying $\Hom(-,A)$ to this resolution shows that $\Hom(K,A)=K(d)$ and $\Ext^1(K,A)$ is 
generated by the vectors $(f x_0) e_i$ with $1 \leq i \leq N$. Applying $\Hom(-,A)$ to 
the short exact sequence $I_D \to I_C \to K$ gives
\[
0 \to \Hom(K,A) \to \Hom(I_C,A) \to \Hom(I_D,A) \to {\Ext}^1(K,A) \to \dots
\]
but $\dim \Hom(K,A) = \dim K(d)=1$, and $\Hom(I_C,A) \cong A(d)$ so that 
$\dim \Hom(I_C,A)={{d+N} \choose {d}}$. Since $\dim {\Ext}^1 (K,A) \leq N$ by the above, 
we conclude that the Hilbert scheme is smooth. The natural rational map 
$\Hilb^{p(z)} (\mathbb P^N) \times  \mathbb P^N \to \Hilb^{p(z)+1} (\mathbb P^N)$ is actually a bijective 
morphism in view of the unique form of the ideal (Remark 
\ref{codim1}), hence an isomorphism by Zariski's main theorem.  
\end{proof}

We are also interested in Hilbert schemes of space curves and obtain the following 
irreducibility result for curves of high genus. Recall that if $C$ is a space curve of degree $d$, 
then $g = p_a (C) \leq {{d-1}\choose{2}}$ with equality for plane curves. 

\begin{thm}\label{irred}
The Hilbert scheme $\Hilb^{dz+1-g} (\mathbb P^3)$ is irreducible for 
$(d,g)$ satisfying 
$d \geq 3$, ${{d-1}\choose{2}}-4 < g \leq {{d-1}\choose{2}}$ and $g > {{d-2}\choose{2}}$. 
\end{thm}

\begin{proof}
The Hilbert scheme is non-empty for all $g \leq {{d-1}\choose{2}}$ due to plane curves 
union isolated points. For $d \geq 3$, the genus of a non-plane curve satisfies 
$g \leq {{d-2}\choose{2}}$ \cite{genus}, so if $C \in \Hilb^{dz+1-g} (\mathbb P^3)$ and 
$C_0 \subset C$ is the curve remaining after removing embedded or isolated points, then $C_0$ is planar, 
hence a complete intersection. Since $C$ is obtained by adding $\leq 3$ embedded or 
isolated points, it is a flat limit of those with isolated points by 
Propositions \ref{onepoint}, \ref{twopoint} and \ref{threepoint} and we conclude 
that the corresponding Hilbert scheme is irreducible. 
\end{proof}

When just one isolated or embedded point is added to a plane curve of 
degree $d$ and genus $g = (d-1)(d-2)/2$, the resulting Hilbert component is 
smooth:

\begin{thm}\label{hd}
The component $H_d \subset \Hilb^{dz+2-g} (\mathbb P^{3})$ of the Hilbert scheme whose 
general member is a degree $d$ plane curve union an isolated point is smooth for all $d \geq 1$. 
Moreover, $H_{d}$ is isomorphic to the blow-up of 
$\Hilb^{dz+1-g} (\mathbb P^{3}) \times \mathbb P^{3}$ along the incidence correspondence.
\end{thm}

\begin{proof}
For $d = 2$ or $3$, this was proved in \cite{CCN} and \cite{PS}, even though $H_d$ is not 
the full Hilbert scheme in these cases. For $d=1$ and $d \geq 4$, $H_d$ is the full Hilbert 
scheme, and it suffices to compute the global sections $H^0({\mathcal N}_D)$ of the normal sheaf 
associated to a point $[D] \in H_d$, so let $D$ be the union of plane curve $C$ and the point 
$p=(0,0,0,1)$. 
If $p \not \in C$, smoothness follows from ${\mathcal N}_D \cong {\mathcal N}_C \oplus {\mathcal N}_p$. 
If $p \in C$ is an embedded point, write $I_C = (z,f)$ with $f \in (x,y)$ and $z=0$ the equation 
of the plane $H$ containing $C$. Consider the exact sequence (\ref{standard}). 
If $\varphi (z) = 0$, then $D \subset H$ and $h^0 ({\mathcal N}_{D,H}) = {d+2 \choose 2}+1$ 
from Theorem \ref{hypersurface}, so the exact sequence 
\[
0 \to {\mathcal N}_{D,H} \to {\mathcal N}_{D,\Pthree} \to \O_D (1) \to 0
\]
yields $h^0({\mathcal N}_{D,\Pthree}) \leq {d+2 \choose 2}+1+h^0(\O_D(1))$.
If $d \geq 4$, then $h^0(\O_D(1))=4$ and we have 
$h^0({\mathcal N}_{D,\Pthree}) \leq \dim H_d$, so $H_d$ is smooth at $D$. 
Similarly, $h^0(\O_D(1))=3$ if $d=1$ and we obtain $h^0({\mathcal N}_D) \leq 7 = \dim H_1$.  

This leaves the possibility that $\varphi(z) \neq 0$. Exact sequence (\ref{standard}) 
shows that $h^1(\I_D(n))=0$ for all $n > 0$, hence the map $(S/I_D)_n \to H^0(\O_D(n))$ 
is an isomorphism for all $n > 0$ ($S$ being the homogeneous coordinate ring of $\Pthree$). 
It follows that the comparison theorem of 
Piene and Schlessinger \cite{PS} applies to $D$ so that $H^0({\mathcal N}_D) \cong \Hom(I_D,S/I_D)_0$. 
Since $\varphi(f)=\lambda \varphi(zw^{d-1})$ for some $\lambda \in k$, 
$\varphi(f-\lambda z w^{d-1})=0$ and we may write $I_D = (xz,yz,z^2,f-\lambda zw^{d-1})$. 
For smoothness at $D$, it suffices to show this when $\lambda=0$, because the members of the 
family parametrized by $\lambda$ are projectively equivalent for $\lambda \neq 0$. 
Thus we may assume $I_D = (xz,yz,z^2,f)$ with $f \in (x,y)$ and write $f = xg+yh$. 

Now consider $\phi \in \Hom(I_D,S/I_D)_0$. Observe that a basis for $(S/I_D)_2$ consists of $\{x^2,xy,y^2,zw,w^2,wx,wy\}$ 
and there is a similar basis for $(S/I_D)_3$ consisting of 11 monomials because $\deg f > 3$. 
In terms of these bases, the Koszul relations $z\phi(xz)=x\phi(z^2),z\phi(xy)=y\phi(z^2),x\phi(yz)=y\phi(xz)$ require that 
\[
\phi(z^2)=a_1 w z, \;\;\;\;\;
\phi(xz) = a_2 wz + a_3 xz + a_4 x^2 + a_5 xy, 
\]
\[
\phi(yz)= a_6 wz + a_3 zy + a_4 xy + a_5 y^2.
\]
Modulo $(xz,yz,z^2)$ we may write 
\[
\phi(f)=a_7 z w^{d-1} + G
\]
with $G \in k[x,y,w]_d$. Now $g \phi(zx) + h \phi(zy) = z \phi(f) = z G$ modulo $I_D$ gives a linear relation 
between the coefficient of $w^d$ in $G$ and $a_2, a_6$. Since $\phi(f)$ is only determined modulo $f$, there are
${d+2 \choose 2}-2$ degrees of freedom in choosing $\phi(f)$, so that 
\[
\dim \Hom(I_D,S/I_D)_0 \leq 7+{d+2 \choose 2} - 2 = {d+2 \choose 2}+5 = \dim H_d.
\]

The second statement follows from Proposition \ref{map} by varying $C$. 
Indeed, the rational map 
$M=\Hilb^{dz+1-g} (\mathbb P^{3}) \times \mathbb P^{3} \to H_{d} \subset \Hilb^{dz+2-g} (\mathbb P^{3})$ 
given by $(C,p) \mapsto C \cup p$ has indeterminacy locus equal to the incidence correspondence 
$\Delta = \{(C,p):p \in C\}$. 
For fixed $C \in \Hilb^{dz+1-g} (\mathbb P^{3})$, the fibre is 
isomorphic to $\mathbb P^{3}$ and via this isomorphism the 
intersection with $\Delta$ is identified with $C \subset \mathbb P^{3}$. 
Thus when $\Delta$ is blown up, the fibre over $C$ is identified with 
${\widetilde {\mathbb P^{3}}} (C)$, which according to Lemma \ref{map} is in bijective 
correspondence with $V \subset H_{d}$ (using the notation in the 
Lemma). It follows that after blowing up the indeterminacy locus $\Delta \subset M$ we obtain a 
bijective map ${\widetilde M}(\Delta) \to H_{d}$, which is an isomorphism by Zariski's main theorem. 
\end{proof}
 
\begin{rmk}\label{2pointsing}{\em
One can verify by similar tangent space calculations that the Hilbert scheme 
of plane curves with two isolated or embedded points is singular exactly along 
the plane curves with double embedded point of type (a) from Proposition \ref{twopoint}. 
It is interesting that the Hilbert scheme is smooth along curves with double embedded 
point of type (b). 
\em}\end{rmk}

\begin{ex}{\em
Theorem \ref{irred} says that $\Hilb^{dz+1-g} (\mathbb P^3)$ is irreducible for  
\[
(d,g) \in \{(3,1),(4,3),(4,2),(5,6),(5,5),(5,4),(6,10),(6,9),(6,8),(6,7)\}.
\]
The pairs $(d,g)=(3,1),(4,3),(5,6)$ and $(6,10)$ correspond to 
plane curves, but the others correspond to plane curves with embedded or isolated points. 
\em}\end{ex}

\begin{ex}{\em
The only locally Cohen-Macaulay curve of degree one is a line. By our results, any curve 
obtained from a line $L$ by adding $\leq 3$ embedded points is a limit of $L$ and the 
right number of isolated points. It follows that $\Hilb^{z+1-g}(\mathbb P^3)$ is irreducible 
of dimension $4-3g$ for $-3 \leq g \leq 0$. On the other hand, it is reducible for 
$g \ll 0$ because the Hilbert scheme of large length schemes is not irreducible \cite{I}.
\em}\end{ex}  
  
\begin{ex}{\em
For curves of degree $2$ and high genus the irreducible components 
of $\Hilb^{2z+1-g} (\mathbb P^3)$ are as follows:

(a) If $g=0$, the Hilbert scheme is irreducible, consisting of plane curves. 

(b) If $g=-1$, there are two irreducible components. The first component $H_1$ 
has general member a pair of skew lines and has dimension 8. The second $H_2$ 
has general member a plane conic union an isolated point and has dimension 11.
There are also plane curves with embedded points, but these lie in $H_2$ 
by Proposition \ref{onepoint}. Both components $H_1$ and $H_2$ are smooth \cite{CCN}.

(c) Similarly if $g=-2$, there are three irreducible components. 
There is the family $H_1$ of double lines of genus $g=-2$ with no embedded points 
of dimension $9$, the family $H_2$ of two skew lines union an isolated point of 
dimension $11$ and the family $H_3$ of conics union two isolated points of dimension 14. 
Because all the underlying locally Cohen-Macaulay curves in question are local complete 
intersections, we know from Remark \ref{onepoint} and Proposition \ref{twopoint} that 
we haven't missed any curves. 

(d) For $g=-3$ we can write down four irreducible components following the same pattern 
as above and our results show that we have missed no irreducible components. However 
when $g=-4$ we cannot be sure that there isn't an irreducible component whose general 
member consists of a plane curves with some horrible quadruple point. 
\em}\end{ex}

\begin{ex}{\em
For curves of degree 3 and high genus, we can make similar lists of the irreducible components 
of $\Hilb^{3z+1-g} (\mathbb P^3)$.

(a) If $g=1$, the Hilbert scheme is irreducible and consists of plane curves. 

(b) If $g=0$, the Hilbert scheme has two irreducible components. The family $H_1$ has 
general member the twisted cubic curve and has dimension 12. The family $H_2$ has general 
member a plane cubic union an isolated point and has dimension 15. This example has been 
well-studied by Piene and Schlessinger \cite{PS}. 

(c) If $g=-1$, there are three irreducible components: the component $H_1$ whose general member 
is a line and a conic has dimension 12; the component $H_2$ whose general member is a twisted 
cubic union an isolated point has dimension 15; the component $H_3$ whose general member is 
a plane cubic union two isolated points has dimension 18. To see that these are all, 
we need to be sure that degenerations of the twisted cubic curve union an embedded point can not 
form an irreducible component of their own, something which is not clear in view of Example 
\ref{onepointfailure} (a). However all ACM curves of degree 3 and 
genus 0 have resolution 
\[
0 \to \O(-3)^2 \to \O(-2)^3 \to \I_C \to 0
\]
\cite{Ell}[Ex. 1] and we can apply Proposition \ref{cilimit}.
\em}\end{ex}

\begin{ex}{\em
Consider the Hilbert schemes $\Hilb^{4z+1-g}(\mathbb P^3)$.

(a) If $g=3$ or $2$, the Hilbert scheme is irreducible by Theorem \ref{irred}. 

(b) If $g=1$, the Hilbert scheme has two irreducible components. One component $H_1$ has general 
member a plane quartic union two isolated points and has dimension 23. Any curves outside of $H_1$ 
have no isolated or embedded points (any non-planar locally Cohen-Macaulay curve satisfies 
the genus bound $g \leq \frac{(d-2)(d-3)}{2}$ \cite{MDP3}), so we are looking at the Hilbert scheme $H_{4,1}$ of locally Cohen-Macaulay curves, which we described in Example \ref{extremal}. 
\em}\end{ex}

This brings us to the last example, which might be known to some
experts, but we have not seen a rigorous proof in the literature. 

\begin{thm}\label{components}
The Hilbert scheme $\Hilb^{4z+1}(\mathbb P^3)$ has 4 irreducible components:
\begin{enumerate}
\item[$H_1$:] The closure of the family of rational quartic curves 
has dimension 16. 

\item[$H_2$:] The family whose general member is a disjoint union 
of a plane cubic and a line has dimension 16.

\item[$H_3$:] The family whose general member is a disjoint union 
of an elliptic quartic curve and a point has dimension 19.

\item[$H_4$:] The family whose general member is a disjoint union 
of a plane quartic curve and 3 distinct points has dimension 26.

\end{enumerate}
\end{thm}

\begin{proof}
The dimension counts are standard, so we only need to show that  
every curve in $\Hilb^{4z+1}(\mathbb P^3)$ is contained in one of these families 
and no family is contained in another. The second part is easy:
family $H_4$ has the largest dimension, but none of the others lie in its closure due to the 
3 isolated/embedded points. Similarly $H_3$ has larger dimension than $H_1,H_2$, but 
$H_1,H_2$ are not in its closure due to the isolated/embedded point. 
Since families $H_1,H_2$ have the same dimension, neither lies in the closure of the other. 

To complete the proof, we show that each curve $[C] \in \Hilb^{4z+1}(\mathbb P^3)$ lies in one of 
the families $H_i$ listed above. Fixing such $C \subset \mathbb P^3$, let 
$C_0 \subset C$ be the purely one-dimensional part of $C$. 
There is no such curve of genus $g=2$ \cite{genus}, leaving three cases. 
If $g(C_0)=0$, then $C=C_0$ is locally Cohen-Macaulay, and 
it is known that the Hilbert scheme $H_{4,0}$ of locally Cohen-Macaulay curves 
has two irreducible components, described in $H_1$ and $H_2$ above \cite{NS}. 
If $g(C_0)=3$, then $C_0$ is a plane quartic and hence a complete intersection. 
It follows from Propositions \ref{onepoint}, \ref{twopoint} and \ref{threepoint} 
that $C$ is a limit of curves which are plane quartics with 3 isolated points, 
so $[C] \in H_4$. If $g(C_0)=1$, then $[C] \in H_3$ by Example \ref{extremal}.
\end{proof}


\begin{thebibliography}{99}

\bibitem{C} D. Chen, Mori's program for the Kontsevich moduli space 
${\overline {\mathcal M}}_{0,0} (\Pthree,3)$, 
Int. Math. Res. Not. IMRN {\bf 2009}, Art. ID rnn 067, 17 pp.

\bibitem{CCN} D. Chen, I. Coskun and S. Nollet, Hilbert scheme  
of a pair of codimension two linear subspaces, to appear in Comm. Algebra.

\bibitem{E} D. Cartwright, D. Erman, M. Velasco and B. Viray, 
Hilbert schemes of 8 points, Algebra Number Theory {\bf 3} (2009) 763--795. 

\bibitem{Ell} G. Ellingsrud, Sur la Schema de Hilbert des Vari\'etes
de Codimension 2 in ${\bf P}^{n}$ a cone de Cohen-Macaulay,
Ann. Scient. \'Ec. Norm. Sup. {\bf 8} (1975) 423-431.

\bibitem{F} J. Fogarty, Algebraic families on an algebraic surface, 
Amer. J. Math. {\bf 90} (1968) 511-Ü521.
 
\bibitem{G} A. Grothendieck, Techniques de construction et th\'eor\`emes d'existence 
en g\'eom\'etrie alg\'ebrique. IV. Les sch\'emas de Hilbert, S\'eminaire Bourbaki, Vol. 6,  
Exp. No. 221, 249--276, Soc. Math. France, Paris, 1995. 

\bibitem{HM} J. Harris and I. Morrison, {\it Moduli of Curves}, 
GTM {\bf 187}, Springer-Verlag, New York, 1998. 

\bibitem{H1} R. Hartshorne, Connectedness of the Hilbert scheme, 
Publ. Math. I.H.E.S. {\bf 29} (1966) 5-48.

\bibitem{AG} R. Hartshorne, {\it Algebraic Geometry} GTM {\bf 52},
Springer-Verlag, Berlin, Heidelberg and New York, 1977.

\bibitem{genus} R. Hartshorne, Genus of Space Curves, Annali Dell' 
Universit\`a di Ferrara, Sezione VII - Scienze Matematiche {\bf 40} 
(1994) 207-223.

\bibitem{Hartconn} R. Hartshorne, On the connectedness of the Hilbert scheme of curves in 
${\mathbb P}\sp 3$, Comm. Algebra  {\bf 28}  (2000) 6059--6077.

\bibitem{I} A. Iarrobino, Reducibility of the family of 0-dimensional 
schemes on a variety, Invent. Math. {\bf 15} (1972) 72--77. 

\bibitem{MDP2} M. Martin-Deschamps and D. Perrin, Sur la Classification 
des Courbes Gauches, Asterisque {184-185}, Soc. Math. de France (1990).

\bibitem{MDP3} M. Martin-Deschamps and D. Perrin, Sur les bornes du module de Rao, 
C. R. Acad. Sci. Paris, t. 137, s\'erie I (1993) 1159--1162.

\bibitem{M} G. Mazzola, Generic finite schemes and Hochschild cocycles,
Comment. Math. Helv. {\bf 55} (1980) 267--293. 

\bibitem{LT} J. Migliore, Introduction to Liaison Theory, 
Progress in Math. 165, Birkh\"auser, 1999.

\bibitem{N} S. Nollet, The Hilbert scheme of degree three curves, 
Ann. Scient. \'Ec. Norm. Sup. 4 serie {\bf 30} (1997) 367--384.

\bibitem{N2} S. Nollet, Deformations of space curves: connectedness 
of Hilbert schemes,  Rend. Semin. Mat. Univ. Politec. Torino  {\bf 64}  (2006) 433--450.

\bibitem{NS} S. Nollet and E. Schlesinger, Hilbert schemes of degree 
four curves, Compositio Math. {\bf 139} (2003) 169--196.

\bibitem{PeS} C. Peskine and L. Szpiro, Liaison des vari\'et\'es 
alg\'ebriques I, Invent. Math. {\bf 26} (1972) 271-302.

\bibitem{PS} R. Piene and M. Schlessinger, 
On the Hilbert scheme compactification of the space of twisted cubics,
Amer. J. Math.  {\bf 107}  (1985) 761--774.

\end{thebibliography}
\end{document}